\documentclass[11pt]{amsart}
\usepackage{amsmath,amsthm}
\newtheorem{theorem}{Theorem}[section]
\newtheorem{cor}[theorem]{Corollary}
\newtheorem{lemma}[theorem]{Lemma}
\newtheorem{prop}[theorem]{Proposition}
\theoremstyle{remark}
\newtheorem{remark}[theorem]{Remark}
\newtheorem{example}[theorem]{Example}
\theoremstyle{definition}
\numberwithin{equation}{section}
\DeclareMathOperator{\Ad}{Ad}
\DeclareMathOperator{\Aut}{Aut}
\DeclareMathOperator{\End}{End}
\DeclareMathOperator{\lsp}{span}
\DeclareMathOperator{\clsp}{\overline{span}}
\newsymbol\rtimes 216F
\newcommand{\rankone}[2]{#1 \otimes\overline{#2}}
\newcommand{\Star}{${}^*$\ndash}

\newcommand{\ndash}{\nobreakdash-}
\newcommand{\field}[1]{\mathbb{#1}}
\newcommand{\CC}{\field{C}}
\newcommand{\NN}{\field{N}}
\newcommand{\TT}{\field{T}}
\newcommand{\ZZ}{\field{Z}}
\newcommand{\Aa}{{\mathcal A}}
\newcommand{\Bb}{{\mathcal B}}
\newcommand{\Ff}{{\mathcal F}}
\newcommand{\Hh}{{\mathcal H}}
\newcommand{\Ii}{{\mathcal I}}
\newcommand{\Kk}{{\mathcal K}}
\newcommand{\Oo}{{\mathcal O}}
\begin{document}
%
% Top Matter
%
\title[Generalized Cuntz algebras]
{Discrete product systems of finite-dimensional Hilbert spaces,
and generalized Cuntz algebras}
\author[Neal J. Fowler]{Neal J. Fowler}
\address{Department of Mathematics  \\
      University of Newcastle\\  NSW  2308\\ AUSTRALIA}
\email{neal@math.newcastle.edu.au}
\date{April 14, 1999}
\thanks{The author thanks the University of Victoria, Canada,
for its hospitality while this research was being conducted.}
\begin{abstract}
To each discrete product system $E$ of finite-dimensional Hilbert spaces
we associate a $C^*$-algebra $\Oo_E$.
When $E$ is the $n$-di\-men\-sion\-al product system over $\NN$,
$\Oo_E$ is the Cuntz algebra $\Oo_n$, and the irrational rotation
algebras appear as $\Oo_E$ for certain one-dimensional product systems
over $\NN^2$.   We give conditions which ensure that $\Oo_E$
is simple, purely infinite, and nuclear.
Our main examples are the lexicographic product systems,
for which we obtain slightly stronger results.
\end{abstract}
\maketitle
%
% Document body
%
\section*{Introduction}

Every representation $\alpha:G\to\Aut B(\Hh)$ of a discrete group $G$
as automorphisms of the algebra $B(\Hh)$
of bounded linear operators on a complex Hilbert space $\Hh$
determines a unique element of $H^2(G,\TT)$:
each $\alpha_s$ is of the form $\Ad U_s$ for some unitary operator $U_s$,
and $U_sU_t = \omega(s,t)U_{st}$ determines a $2$\ndash cocycle
$\omega$ whose cohomology class is independent of the choice
of implementing unitaries.
If instead one starts with an action
of a semigroup $P$ as endomorphisms of $B(\Hh)$,
then, as described in \cite[Remark~2.3]{fowrae},
the cohomological obstruction one encounters is a
{\em product system\/} over $P$:
roughly speaking,  a collection
$E = \{E_s:s\in P\}$ of complex Hilbert spaces together with
an associative multiplication which implements unitary isomorphisms
$E_s\otimes E_t \to E_{st}$.
Product systems were first defined by Arveson \cite{arv}
in his study of one-parameter semigroups
of endomorphisms of $B(\Hh)$.

In this note we continue the investigations of
\cite{dinhjfa, dinhjot, laca, fowrae, fowler}
into the $C^*$\ndash algebras
associated with discrete product systems.
Each of the algebras studied in these previous papers
can be regarded as a Toeplitz algebra:
unit vectors of the product system
correspond to isometries in the universal $C^*$\ndash algebra,
and finite orthonormal sets of vectors from a fiber $E_s$
map to isometries whose range projections are orthogonal
and whose sum is strictly less than the identity.
Here we consider only product systems $E$ whose fibers
are finite-dimensional,
and in our universal $C^*$\ndash algebra $\Oo_E$
the range projections of the isometries associated
with an orthonormal basis for any fiber sum to the identity.
Hence $\Oo_E$ is generated by a collection
of Cuntz algebras, one for each element of the underlying semigroup $P$;
the relations between these Cuntz algebras
are determined by the multiplication in $E$.
We think of $\Oo_E$ as a semigroup version
of the twisted group algebra $C^*(G,\omega)$.

We begin in Section~\ref{section:OE} by establishing
the existence of $\Oo_E$ and discussing a few examples.
In particular one can recover as $\Oo_E$ the Cuntz algebras
(Remark~\ref{remark:cuntz algebra}) and the irrational rotation
algebras (Remark~\ref{remark:irrational}).
In Section~\ref{section:abelian} we develop some basic results
regarding the structure of $\Oo_E$
when the underlying semigroup $P$ embeds in an abelian group;
this is a key hypothesis in most of our results.
In Section~\ref{section:simplicity} we prove our main result,
Theorem~\ref{theorem:main}:
if no two fibers of $E$ have the same dimension,
then $\Oo_E$ is simple and purely infinite.
For the lexicographic product systems of \cite{fowrae}
our Theorem~\ref{theorem:lexicographic} gives a sharper simplicity  result:
$\Oo_E$ is simple  if and only if the dimension function is injective.
This shows in particular that $\Oo_E$ need not be simple
even if the only one-dimensional fiber of $E$ is the fiber over
the identity of $P$.
For product systems over $\NN^2$ we show in Theorem~\ref{theorem:N2}
that $\Oo_E$ is either
simple or isomorphic to $\Oo_l\otimes C(\TT)$, where $l$ is determined
by the dimensions of the fibers over $(1,0)$ and $(0,1)$.
Finally, in Theorem~\ref{theorem:nuclearity} we generalize a technique
of Laca \cite{laca} to show that $\Oo_E$ is nuclear
if $E$ admits a twisted unit.

\section{The Cuntz algebra of a product system}\label{section:OE}

Suppose $P$ is a countable discrete semigroup with identity $e$.
A {\em product system\/} over $P$ is a family $p:E\to P$
of nontrivial separable complex Hilbert spaces $E_t:= p^{-1}(t)$
which is endowed with an associative multiplication $E\times E\to E$
in such a way that $p$ is a semigroup homomorphism,
and such that for every $s,t\in P$ the map
$x\otimes y\in E_s\otimes E_t \mapsto xy\in E_{st}$ 
extends to a unitary isomorphism.
We also insist that $\dim E_e = 1$,
so that $E$ has an identity $1\in E_e$ \cite[Lemma~1.3]{fowrae}.

A {\em representation\/} of $E$ in a $C^*$\ndash algebra $B$
is a map $\phi:E\to B$ which satisfies
\begin{itemize}
\item[(i)] $\phi(xy) = \phi(x)\phi(y)$ for every $x,y\in E$; and
\item[(ii)] $\phi(y)^*\phi(x) = (x\mid y) \phi(1)$
if $p(x) = p(y)$.
\end{itemize}
This definition is slightly different than that found in
\cite{fowrae}, where representations were on Hilbert space
and condition (ii) was that
$\phi(y)^*\phi(x) = (x\mid y) 1$ if $p(x) = p(y)$.
The main advantage to our definition is that it allows us to
consider the trivial map $E\to\{0\}$ as a representation,
without regarding $\{0\}$ as a unital $C^*$\ndash algebra.
If $\phi \ne 0$, then $\phi(1)$ is a nontrivial projection
which serves as an identity for the $C^*$\ndash subalgebra
$C^*(\phi)\subseteq B$ generated by $\phi(E)$;
moreover, $\phi$ restricts to
an isometric linear map on each of the fibers of $E$
\cite[p.8]{arv}.

For nontrivial $\phi$,
condition (ii) implies that each $\phi(x)$ is a multiple of an isometry
in $C^*(\phi)$;
indeed, if $\Bb_t$ is an orthonormal basis for $E_t$,
then $\{\phi(f): f\in \Bb_t\}$ is a family of isometries with
mutually orthogonal range projections.
In this note we will consider only product systems whose fibers
are finite-dimensional, and we are primarily interested in representations
$\phi:E\to B$ which satisfy
\begin{equation}\label{eq:cuntz rep}
\sum_{f\in \Bb_t} \phi(f)\phi(f)^* = \phi(1)
\qquad\text{for all $t\in P$.}
\end{equation}
We call a representation $\phi$ which satisfies \eqref{eq:cuntz rep}
a {\em Cuntz representation\/}.
While nontrivial representations of $E$ always exist 
\cite[Lemma~1.10]{fowrae},
it is not clear that every product system admits a nontrivial
Cuntz representation.

\begin{prop}\label{prop:OE}
Let $E$ be a product system over $P$ of finite-dimensional Hilbert spaces.
There is a pair $(\Oo_E,i_E)$ consisting of a $C^*$\ndash algebra $\Oo_E$
and a Cuntz representation $i_E:E\to\Oo_E$ with the following properties:

\textup{(a)} for every Cuntz representation $\phi$ of $E$,
there is a homomorphism $\phi_*$ of $\Oo_E$,
called the {\em integrated form of $\phi$\/},
such that $\phi_*\circ i_E = \phi$; and

\textup{(b)} $\Oo_E$ is generated as a $C^*$\ndash algebra by $i_E(E)$.

\noindent The pair $(\Oo_E,i_E)$ is unique up to canonical isomorphism.
\end{prop}

\begin{remark}\label{remark:trivial}
We emphasize that $\Oo_E$ is trivial if $E$ does not admit
a nontrivial Cuntz representation.
In all other cases $\Oo_E$ is unital and $i_E$ is isometric.
\end{remark}

\begin{proof}[Proof of Proposition~\ref{prop:OE}]
The Proposition can be proved by modifying the standard argument
of, for example, \cite[Proposition~1.3]{fowrae2}.
Briefly, let $S$ be a collection of Cuntz representations on
Hilbert space which are cyclic (i.e. generate a $C^*$\ndash algebra
which admits a cyclic vector), and such that every cyclic
Cuntz representation is unitarily equivalent to an element of $S$.
Define $i_E := \bigoplus_{\phi\in S} \phi$ and $\Oo_E := C^*(i_E)$.
Condition (a) holds because every Cuntz representation of $E$
on a Hilbert space decomposes as the direct sum of a zero representation
and a collection of cyclic representations,
and the uniqueness assertion follows from a standard argument.
\end{proof}

\begin{example}\label{example:lexicographic}
{\em (Lexicographic product systems.)\/}
For any product system $E$ whose fibers are finite-dimensional,
$s \mapsto \dim E_s$ is a semigroup homomorphism from $P$
to the multiplicative positive integers $\NN^*$.
We call this homomorphism the {\em dimension function\/} of $E$.
Let $d:P\to\NN^*$ be an arbitrary homomorphism.
In \cite[Examples~1.4(E2)]{fowrae} it was shown how one can construct
a product system $E(d)$, called the {\em lexicographic product system
determined by $d$\/}, whose dimension function is $d$:
for each $n$ let $\{\delta_0,\dots,\delta_{n-1}\}$
be the canonical basis for  $\CC^n$,
take
\[
E(d) := \bigsqcup_{s\in P} \{s\}\times \CC^{d(s)},
\]
and define multiplication on basis vectors using
the lexicographic order on
$\{0,\dots, d(r) - 1\} \times \{0, \dots, d(s) - 1\}$; that is,
\[
(r,\delta_j)(s,\delta_k) := (rs,\delta_{jd(s) + k})
\qquad\text{for $0 \le j \le d(r) - 1$ and $0 \le k \le d(s) - 1$.}
\]
Multiplication in $E(d)$ is defined by extending bilinearly.

There is a distinguished Cuntz representation of $E(d)$
on $L^2(\TT)$.  Suppose $r\in P$ and $0 \le j \le d(r) - 1$.
Let $S(r,\delta_j)$ be the operator on $L^2(\TT)$
whose value on a vector $\xi\in L^2(\TT)$ is given by
\[
S(r,\delta_j)\xi(e^{2\pi it})
:=  \begin{cases}
    d(r)^{1/2} \xi(e^{2\pi i(td(r) - j)})
      & \text{if $t\in \bigl[ \frac j{d(r)}, \frac{j+1}{d(r)} \bigr)$} \\
    0 & \text{otherwise.} \end{cases}
\]
It is easy to see that $S(r,\delta_j)$ is an isometry,
and that
\[
\sum_{j=0}^{d(r) - 1} S(r,\delta_j)S(r,\delta_j)^* = 1 = S(1).
\]
Since
\begin{align*}
& S(r,\delta_j)S(s,\delta_k)\xi(e^{2\pi it}) \\
& \qquad = \begin{cases}
    d(r)^{1/2} S(s,\delta_k)\xi(e^{2\pi i(td(r) - j)})
      & \text{if $t\in \bigl[ \frac j{d(r)}, \frac{j+1}{d(r)} \bigr)$} \\
    0 & \text{otherwise} \end{cases} \\
& \qquad = \begin{cases}
    d(r)^{1/2} d(s)^{1/2} \xi(e^{2\pi i((td(r) - j)d(s) - k)})
      & \text{if $td(r) - j\in \bigl[ \frac k{d(s)}, \frac{k+1}{d(s)} \bigr)$} \\
    0 & \text{otherwise} \end{cases} \\
& \qquad = \begin{cases}
    d(rs)^{1/2} \xi(e^{2\pi i(td(rs) - (jd(s) + k))})
      & \text{if $t \in \bigl[
              \frac {jd(s) + k}{d(rs)}, \frac{jd(s) + k + 1}{d(rs)} \bigr)$} \\
    0 & \text{otherwise} \end{cases} \\
& \qquad = S(rs,\delta_{jd(s)+k})\xi(e^{2\pi it}),
\end{align*}
we have
$S(r,\delta_j)S(s,\delta_k) = S(rs,\delta_{jd(s)+k})
= S((r,\delta_j)(s,\delta_k))$.
Hence defining
\[
S(r,x) := \sum_{j=0}^{d(r)-1} (x \mid \delta_j) S(r,\delta_j)
\qquad\text{for $(r,x)\in E(d)$}
\]
gives a Cuntz representation $S:E(d)\to B(L^2(\TT))$.
\end{example}

\begin{remark}\label{remark:cuntz algebra}
Fix $n \ge 2$ and let $d$ be the homomorphism $a\in\NN \mapsto n^a\in\NN^*$.
Then $\Oo_{E(d)}$ is the Cuntz algebra $\Oo_n$.
\end{remark}

\begin{example}
{\em (One-dimensional product systems.)\/}
Suppose $\omega$ is a multiplier of $P$;
that is, a function  $\omega:P\times P\to\TT$
such that $\omega(e,e) = 1$ and
\[
\omega(r,s)\omega(rs,t) = \omega(r,st)\omega(s,t)
\qquad\text{for all $r,s,t\in P$.}
\]
Then
$(r,z)(s,w) := (rs,\omega(r,s)zw)$
for $r,s\in P$ and $z,w\in\CC$
defines an associative multiplication which gives
$P\times\CC$ the structure of a product system;
we write $(P\times\CC)^\omega$ for this product system.
It is easy to see that every product system whose fibers
are one-dimensional is isomorphic to $(P\times\CC)^\omega$
for some multiplier $\omega$.

Suppose $P$ is an Ore semigroup; that is,
suppose $P$ is cancellative and satisfies $Pr\cap Ps \ne\emptyset$
for every pair $r,s\in P$.  (For example, every cancellative commutative
semigroup has this property.)
By \cite[Theorem~1.1.2]{laca2}, $P$ can be embedded
in a group $G$ with $P^{-1}P = G$,
and by \cite[Theorem~1.2.2]{laca2} there is a multiplier
$\omega'$ of $G$ which extends $\omega$.
Define $l^\omega:(P\times\CC)^\omega\to B(\ell^2(G))$
by
\[
l^\omega(r,z)\xi(s) := z\omega'(r,r^{-1}s)\xi(r^{-1}s)
\qquad\text{for $\xi\in\ell^2(G)$ and $s\in G$.}
\]
It is routine to check that each $l^\omega(r,1)$ is unitary
and that $l^\omega$ is multiplicative.
Hence $l^\omega$ is a Cuntz representation of $(P\times\CC)^\omega$.
\end{example}

\begin{remark}\label{remark:irrational}
Let $\theta\in(0,1)$ be irrational,
let  $\omega:\NN^2\times\NN^2\to\TT$ be the multiplier
$\omega((a,b),(c,d)) := e^{2\pi i\theta bc}$,
and let $E = (\NN^2\times\CC)^\omega$.
Then $U := i_E((0,1),1)$ and $V := i_E((1,0),1)$
are unitaries which generate $\Oo_E$ and satisfy
$UV = e^{2\pi i\theta} VU$,
so $\Oo_E$ is the irrational rotatation algebra $\Aa_\theta$.
\end{remark}

\begin{example}
{\em (Twisting.)\/} If $E$ is a product system over $P$
and $\omega$ is a multiplier of $P$,
then $(x,y)\in E\times E \mapsto \omega(p(x),p(y))xy$
defines a multiplication on $E$
which also gives $E$ the structure of a product system;
we write $E^\omega$ for this new system,
and say $E$ has been {\em twisted by $\omega$\/}.
If $\phi$ is a Cuntz representation of $E$
on a Hilbert space $\Hh$, then
\[
\phi^\omega(x) := \phi(x)\otimes l^\omega(p(x),1)
\qquad\text{for $x\in E^\omega$}
\]
defines a Cuntz representation $\phi^\omega$
of $E^\omega$ on $\Hh\otimes\ell^2(G)$.
\end{example}

The results obtained in the previous examples allow us to state:

\begin{prop}\label{prop:nontrivial}
Every twisted lexicographic product system over an Ore semigroup
admits a nontrivial Cuntz representation.
\end{prop}

\section{Product systems over abelian semigroups}
\label{section:abelian}

The following Proposition collects some results concerning the
structure of $\Oo_E$ when $P$ embeds in an abelian group.

\begin{prop}\label{prop:abelian}
Suppose $P$ is a subsemigroup of a countable abelian group $G$
and $E$ is a product system over $P$ of finite-dimensional Hilbert spaces.

\textup{(1)}
Let $\phi$ be a Cuntz representation of $E$,
let $s,t\in P$, and let $\Bb_s$ and $\Bb_t$
be orthonormal bases for $E_s$ and $E_t$, respectively.
Then for any $x'\in E_s$ and $y'\in E_t$ we have
\[
\phi(y')^*\phi(x') = \sum_{x\in \Bb_s} \sum_{y\in \Bb_t} (x'y \mid y'x) \phi(x)\phi(y)^*.
\]

\textup{(2)} $\Oo_E = \clsp\{i_E(x)i_E(y)^*: x,y\in E\}$.

\textup{(3)} There is a strongly continuous action $\gamma:\widehat G\to\Aut\Oo_E$,
called the {\em gauge action\/},
such that
\[
\gamma_\lambda(i_E(x)) = \lambda(p(x))i_E(x)
\qquad\text{for all $\lambda\in\widehat G$ and $x\in E$.}
\]

\textup{(4)} Let $m$ be Haar measure on $\widehat G$.
Then
\[
\Phi(b) := \int_{\widehat G} \gamma_\lambda(b)\,dm(\lambda)
\qquad\text{for $b\in\Oo_E$}
\]
defines a faithful conditional expectation $\Phi$
of $\Oo_E$ onto $\Oo_E^\gamma$,
the fixed-point algebra of the gauge action.

\textup{(5)} $\Oo_E^\gamma
= \clsp\{i_E(x)i_E(y)^*: x,y\in E,\ p(x) = p(y)\}$.
\end{prop}

\begin{proof}
(1) We use \eqref{eq:cuntz rep} to calculate
\begin{align*}
\phi(y')^*\phi(x')
& = \Bigl( \sum_{x\in \Bb_s} \phi(x)\phi(x)^* \Bigr)
    \phi(y')^*\phi(x')
    \Bigl( \sum_{y\in \Bb_t} \phi(y)\phi(y)^* \Bigr) \\
& = \sum_{x\in \Bb_s} \sum_{y\in \Bb_t}
    \phi(x)\phi(y'x)^*\phi(x'y)\phi(y)^* \\
& = \sum_{x\in \Bb_s} \sum_{y\in \Bb_t}
    (x'y \mid y'x)\phi(x)\phi(y)^*,
\end{align*}
noting that $p(y'x) = p(x'y)$ since $P$ is abelian.

(2) Take $\phi = i_E$ in (1) to see that
the linear span of monomials of the form $i_E(x)i_E(y)^*$
is closed under multiplication, and hence a \Star algebra.

(3) For each $\lambda\in\widehat G$
the map $x\in E \mapsto \lambda(p(x))i_E(x)$
is a Cuntz representation, and hence integrates to
an endomorphism $\gamma_\lambda$ of $\Oo_E$.
Since $\gamma_\lambda \circ \gamma_{\lambda^{-1}}$
and $\gamma_{\lambda^{-1}} \circ \gamma_\lambda$
are both the identity on $\Oo_E$, $\gamma_\lambda$
is an automorphism.  Obviously $\gamma$ is a group
homomorphism, and its continuity follows from a straightforward
$\epsilon/3$ argument.

(4) This is a standard result about automorphic actions of
discrete abelian groups.

(5) Since $\Phi(i_E(x)i_E(y)^*) = \delta_{p(x),p(y)}i_E(x)i_E(y)^*$,
(5) follows from (2) and the continuity of $\Phi$.
\end{proof}

Our next goal is to give an abstract characterization of the
fixed-point algebra $\Oo_E^\gamma$.  First some notation.
If $s,t\in P$, $S\in\Kk(E_s)$ and $T\in\Kk(E_t)$,
write $S\otimes T$ for the operator on $E_{st}$ which satisfies
\[
S\otimes T(xy) = (Sx)(Ty)
\qquad\text{for $x\in E_s$ and $y\in E_t$.}
\]
Write $1^t$ for the identity operator on $E_t$.

Define a relation $\preceq$ on $P$ by $s \preceq t$ if and only if $t\in sP$,
and observe that $\preceq$ is a preorder on $P$;
that is, it is reflexive and transitive.
Since $P$ is commutative, for every pair $s,t\in P$
we have $st \in sP\cap tP$, and hence $(P,\preceq)$ is upwardly directed.
Now $S\mapsto S\otimes 1^t$
is a unital embedding of $\Kk(E_s)$ in $\Kk(E_{st})$,
and since multiplication in $E$ is associative we have
$S\otimes 1^{tr} = (S\otimes 1^t)\otimes 1^r$.
Hence $S\in\Kk(E_s) \mapsto S\otimes 1^t\in\Kk(E_{st})$
is a directed system of $C^*$\ndash algebras,
and we can define
\[
\Ff_E := \varinjlim \Kk(E_s).
\]
Since each $\Kk(E_s)$ is simple, so is $\Ff_E$.
(It is not hard to see that $\Ff_E$ is a UHF algebra.)
Let $\iota_s$ be the canonical embedding of $\Kk(E_s)$ in $\Ff_E$.

\begin{prop}\label{prop:FE}
There is a unique homomorphism $i_\Ff:\Ff_E\to\Oo_E$
which satisfies
\[
i_\Ff(\iota_{p(x)}(\rankone xy)) = i_E(x)i_E(y)^*
\qquad\text{whenever $p(x) = p(y)$,}
\]
where $\rankone xy$ is the rank-one operator
$z\in E_{p(x)} \mapsto (z\mid y)x$.
The image of $\Ff_E$ is precisely $\Oo_E^\gamma$.
If $E$ admits a nontrivial Cuntz representation,
then $i_\Ff$ is injective.
\end{prop}

\begin{proof}
For each $s\in P$, there is a unique homomorphism $\sigma_s:\Kk(E_s)\to\Oo_E$
such that $\sigma_s(\rankone xy) = i_E(x)i_E(y)^*$ for all $x,y\in E_s$.
Let $\Bb_t$ be an orthonormal basis for $E_t$.  Since
\begin{align*}
\sigma_{st}((\rankone xy)\otimes  1^t)
& = \sigma_{st}\Bigl( \sum_{f\in \Bb_t} (\rankone xy) \otimes (\rankone ff) \Bigr)
  = \sigma_{st}\Bigl( \sum_{f\in \Bb_t} (\rankone {xf}{yf} \Bigr) \\
& = \sum_{f\in \Bb_t} i_E(xf)i_E(yf)^*
  = i_E(x)\Bigl( \sum_{f\in \Bb_t} i_E(f)i_E(f)^* \Bigr) i_E(y)^* \\
& = i_E(x)i_E(y)^* = \sigma_s(\rankone xy),
\end{align*}
we deduce that $\sigma_{st}(S\otimes 1^t) = \sigma_s(S)$ for all $S\in\Kk(E_s)$,
and hence there is a homomorphism $i_\Ff:\Ff_E\to\Oo_E$
such that $i_\Ff\circ\iota_s = \sigma_s$ for every $s\in P$.
From Proposition~\ref{prop:abelian}(5) it is obvious that
$i_\Ff$ maps $\Ff_E$ onto $\Oo_E^\gamma$.
If $E$ admits a nontrivial Cuntz representation then $i_E$ is nonzero
(see Remark~\ref{remark:trivial}), hence each $\sigma_s$
is nonzero, and finally $i_\Ff$ is nonzero.
Since $\Ff_E$ is simple, we deduce that $i_\Ff$ is injective.
\end{proof}

\section{Simplicity and pure infiniteness}
\label{section:simplicity}

\begin{theorem}\label{theorem:main}
Let $G$ be a countable abelian group,
let $P$ be a subsemigroup of $G$ which contains the identity,
and let $E$ be a product system over $P$ of finite-dimensional Hilbert spaces
which admits a nontrivial Cuntz representation.
If the dimension function $s \mapsto \dim E_s$ is injective,
then $\Oo_E$ is simple and purely infinite.
\end{theorem}

To prove the Theorem we require a technical lemma.
For the Proposition which precedes it, see
\cite[Proposition~2.7]{arv} and \cite[Proposition~1.11 and Lemma~3.6]{fowrae}.

\begin{prop}\label{prop:alpha}
Let $E$ be a product system over $P$ of finite-dimensional Hilbert spaces,
and let $\phi$ be a Cuntz representation of $E$
in a unital $C^*$\ndash algebra $B$ such that $\phi(1) = 1$.
For each $s\in P$ there is a unital endomorphism $\alpha^\phi_s$ of $B$
which satisfies
\[
\alpha^\phi_s(b) = \sum_{f\in \Bb_s} \phi(f)b\phi(f)^*
\qquad\text{for all $b\in B$}
\]
whenever $\Bb_s$ is an orthonormal basis for $E_s$.
Moreover, $\alpha_e$ is the identity, and
\[
\alpha^\phi_{st}(b)\phi(x) = \phi(x)\alpha^\phi_t(b)
\qquad\text{for all $x\in E_s$.}
\]
\end{prop}

\begin{lemma}\label{lemma:kill}
Let $E$ be a product system satisfying the hypotheses
of Theorem~\ref{theorem:main},
and let $\phi$ be a Cuntz representation of $E$
in a unital $C^*$\ndash algebra $B$ such that $\phi(1) = 1$.
Suppose $x_1$, \dots, $x_n$,
$y_1$, \dots, $y_n\in E$ satisfy $p(x_i) \ne p(y_i)$ for $1 \le i \le n$.
Let $c\in P$ be such that $p(x_i)^{-1}c \in P$ and $p(y_i)^{-1}c\in P$
for every $i$.
Then there is a unit vector $w\in E$ such that
\begin{equation}\label{eq:w}
\alpha^\phi_c(\phi(w)\phi(w)^*)\phi(x_i)\phi(y_i)^*
\alpha^\phi_c(\phi(w)\phi(w)^*) = 0
\qquad\text{for $1 \le i \le n$.}
\end{equation}
\end{lemma}

\begin{proof}
Define $s_i := p(x_i)^{-1}c$, $t_i := p(y_i)^{-1}c$,
and for each $r\in P$ let $\Bb_r$ be an orthonormal basis for $E_r$.
Suppose $Q\in \phi_*(\Oo_E)$;
we will eventually take $Q = \phi(w)\phi(w)^*$.
By Proposition~\ref{prop:alpha} we have
\begin{multline*}
\alpha^\phi_c(Q)\phi(x_i)\phi(y_i)^*\alpha^\phi_c(Q)
 = \phi(x_i)\alpha^\phi_{s_i}(Q)\alpha^\phi_{t_i}(Q)\phi(y_i)^* \\
 = \sum_{f'\in \Bb_{s_i}} \sum_{g'\in \Bb_{t_i}}
    \phi(x_i)\phi(f')Q\phi(f')^*\phi(g')Q\phi(g')^*\phi(y_i)^*,
\end{multline*}
and by Proposition~\ref{prop:abelian}(1)
\[
Q\phi(f')^*\phi(g')Q
 = \sum_{g\in \Bb_{t_i}} \sum_{f\in \Bb_{s_i}}
    (g'f \mid f'g) Q\phi(g)\phi(f)^*Q,
\]
so it suffices to find a unit vector $w\in E$ such that
\begin{equation}\label{eq:sufficient}
\phi(w)^*\phi(g)\phi(f)^*\phi(w) = 0
\qquad\text{whenever $(f,g)\in \Bb_{s_i}\times \Bb_{t_i}$
for some $i$.}
\end{equation}
Let $(f_l, g_l)_{l=1}^m$ be an enumeration of these pairs.
We claim that for each $j\in \{0, \dots, m\}$
there is a unit vector $v_j \in E$ such that
\begin{equation}\label{eq:v}
\phi(v_j)^*\phi(g_l)\phi(f_l)^*\phi(v_j) = 0
\qquad\text{for $1 \le l \le j$.}
\end{equation}
Given the claim, $w := v_m$ satisfies \eqref{eq:sufficient},
and hence \eqref{eq:w}, completing the proof. 
The claim is vacuous when $j = 0$: taking any unit vector $v_0\in E$ works.
Suppose inductively that
there exists $k \le m-1$ such that \eqref{eq:v} holds when $j = k$.
Let  $r := p(v_k)$, $s := p(f_{k+1})$, and $t := p(g_{k+1})$. 
Since $s \ne t$, by hypothesis $\dim E_s \ne \dim E_t$,
and we can assume without loss of generality that $\dim E_s < \dim E_t$.
We have
\begin{align*}
& \phi(v_k)^*\phi(g_{k+1})\phi(f_{k+1})^*\phi(v_k) \\
& \quad = \sum_{f\in \Bb_s} \sum_{g\in \Bb_t} \sum_{v\in \Bb_r}
           \phi(g)\phi(g)^*\phi(v_k)^*\phi(g_{k+1})
           \phi(v)\phi(v)^*\phi(f_{k+1})^*\phi(v_k)
           \phi(f)\phi(f)^* \\
& \quad = \sum_{f\in \Bb_s} \sum_{g\in \Bb_t} \sum_{v\in \Bb_r}
           (g_{k+1}v \mid v_kg)(v_kf \mid f_{k+1}v) \phi(g)\phi(f)^* \\
& \quad = \sum_{f\in \Bb_s} \phi(u_f)\phi(f)^*,
\end{align*}
where
\[
u_f := \sum_{g\in \Bb_t} \sum_{v\in \Bb_r}
           (g_{k+1}v \mid v_kg)(v_kf \mid f_{k+1}v) g \in E_t.
\]
Since $\dim E_s < \dim E_t$, there is a unit vector
$v'\in E_t$ which is orthogonal to each $u_f$,
and taking $v_{k+1} := v_kv'$ gives \eqref{eq:v} for $j = k+1$.
\end{proof}

\begin{proof}[Proof of Theorem~\ref{theorem:main}]
Suppose $B$ is a $C^*$\ndash algebra and
$\pi:\Oo_E\to B$ is a nonzero homomorphism.
We will show that $\pi$ is injective,
thus establishing the simplicity of $\Oo_E$.
It does not harm to assume that $\pi$ is surjective,
so that $B$ is unital and $\pi(1) = 1$.

Recall from Proposition~\ref{prop:abelian}(4) that there
is a faithful expectation $\Phi$ of $\Oo_E$ onto $\Oo_E^\gamma$,
given by averaging over the orbits of the gauge action.
We will show that there is an expectation $\Phi_\pi$
of $\pi(\Oo_E)$ onto $\pi(\Oo_E^\gamma)$ such that
\begin{equation}\label{eq:spatial Phi}
\Phi_\pi\circ\pi = \pi\circ\Phi.
\end{equation}
To see that this implies that $\pi$ is injective,
suppose $b\in\ker\pi$.
Then
$\pi\circ\Phi(b^*b) = \Phi_\pi\circ\pi(b^*b) = 0$,
so $\Phi(b^*b)\in \Oo_E^\gamma\cap\ker\pi$.
But $\Oo_E^\gamma$ is simple and contains the identity of $\Oo_E$,
so $\Oo_E^\gamma\cap\ker\pi = \{0\}$.
Thus $\Phi(b^*b) = 0$, and we deduce from the faithfulness
of $\Phi$ that $b = 0$.

By Proposition~\ref{prop:abelian}(2),
finite sums of the form $b = \sum i_E(x_i)i_E(y_i)^*$
are dense in $\Oo_E$.
Hence to prove the existence of an expectation $\Phi_\pi$
satisfying \eqref{eq:spatial Phi}, it suffices to fix such an
element $b$ and show that
$\lVert \pi(b) \rVert \ge \lVert \pi(\Phi(b)) \rVert$.
Thus with $\phi := \pi\circ i_E$, we must show that
\begin{equation}\label{eq:inequality}
\lVert \sum \phi(x_i)\phi(y_i)^* \rVert
\ge
\Bigl\lVert \sum_{p(x_i) = p(y_i)} \phi(x_i)\phi(y_i)^* \Bigr\rVert.
\end{equation}

Let $c := \prod_i p(x_i)p(y_i)$.
Since $\phi$ is a Cuntz representation of $E$,
the relations \eqref{eq:cuntz rep}
allow us to assume that $p(x_i) = p(y_i) = c$
whenever $p(x_i) = p(y_i)$.
Since $\phi(1) = 1$,
Lemma~\ref{lemma:kill} applies
and provides a unit vector $w\in E$ such that
\[
\alpha^\phi_c(\phi(w)\phi(w)^*)
\phi(x_i)\phi(y_i)
\alpha^\phi_c(\phi(w)\phi(w)^*)
= 0
\qquad\text{whenever $p(x_i) \ne p(y_i)$.}
\]
Let $Q := \phi(w)\phi(w)^*$.
Then
\begin{equation}\label{eq:Q inequality}
\begin{split}
\lVert \sum \phi(x_i)\phi(y_i)^* \rVert
& \ge  \lVert \alpha_c(Q) \sum \phi(x_i)\phi(y_i)^* \alpha_c(Q) \rVert \\
& = \Bigl\lVert \alpha_c(Q)\sum_{p(x_i) = p(y_i)}
    \phi(x_i)\phi(y_i)^*\alpha_c(Q)
    \Bigr\rVert \\
& = \Bigl\lVert \sum_{p(x_i) = p(y_i)} \phi(x_i)Q\phi(y_i)^*
    \Bigr\rVert.
\end{split}
\end{equation}
Since $\rankone xy \mapsto \phi(x)Q\phi(y)^*$
and $\rankone xy \mapsto \phi(x)\phi(y)^*$
both extend linearly to nontrivial homomorphisms
of the simple algebra $\Kk(E_c)$,
we have
\[
\Bigl\lVert \sum_{p(x_i) = p(y_i)} \phi(x_i)Q\phi(y_i)^* \Bigr\rVert
= \Bigl\lVert \sum_{p(x_i) = p(y_i)} \phi(x_i)\phi(y_i)^* \Bigr\rVert.
\]
Combining this with \eqref{eq:Q inequality}
gives \eqref{eq:inequality},
completing the proof of simplicity.

Our proof that $\Oo_E$ is purely infinite is an adaptation of
the proof of \cite[Proposition~5.3]{bprs}.
Let $A$ be a hereditary subalgebra of $\Oo_E$; we will show that
$A$ has an infinite projection.
Fix a positive element $a\in A$, scaled so that $\lVert\Phi(a)\rVert = 1$.
Choose a finite sum
$b = \sum i_E(x_i)i_E(y_i)^* \in\Oo_E$
such that $b \ge 0$ and $\lVert a - b \rVert < 1/4$.
Then $b_0 := \Phi(b)$ is also positive and satisfies $\lVert b_0 \rVert \ge 3/4$.
By applying the relations \eqref{eq:cuntz rep}
we can assume that there  exists $c\in P$
such that $p(x_i) = p(y_i) = c$ whenever $p(x_i) = p(y_i)$, and such that
$p(x_i)^{-1}c \in P$ and $p(y_i)^{-1}c \in P$ for all $i$.
Then $b_0$ is a positive element of the algebra
\[
\Ff_c := \lsp\{i_E(x)i_E(y)^*: p(x) = p(y) = c\},
\]
and hence its image under the canonical isomorphism $\Ff_c \cong \Kk(E_c)$
has a unit eigenvector $f\in E_c$ with eigenvalue $\lVert b_0 \rVert$.
It follows that $b_0i_E(f) = \lVert b_0 \rVert i_E(f)$,
and hence the projection $r := i_E(f)i_E(f)^*$ satisfies
$rb_0r = \lVert b_0 \rVert r$.

By Lemma~\ref{lemma:kill}, there is a unit vector $w\in E$ such that
\[
\alpha^{i_E}_c(i_E(w)i_E(w)^*)
i_E(x_i)i_E(y_i)
\alpha^{i_E}_c(i_E(w)i_E(w)^*)
= 0
\]
whenever $p(x_i) \ne p(y_i)$.
Let $q := i_E(fw)i_E(fw)^*$,
noting that $q \le r$ and $q \le \alpha^{i_E}_c(i_E(w)i_E(w)^*)$.
Then
\[
qbq = qb_0q = qrb_0rq
= \lVert b_0 \rVert qrq
= \lVert b_0 \rVert q
\ge \textstyle\frac 34 q,
\]
and since $\lVert a - b \rVert < 1/4$ we thus have
$qaq \ge qbq - \frac 14q \ge \frac 12q$.
It follows that $qaq$ is invertible in $q\Oo_Eq$.  Let $c$ be its inverse,
and define $v := c^{1/2}qa^{1/2}$.
Then $v$ is a partial isometry since $vv^* = c^{1/2}qaqc^{1/2} = q$,
and its initial projection belongs to $A$ since
$v^*v = a^{1/2}qcqa^{1/2} \le \lVert c\rVert a$
and $A$ is hereditary.
Since $v^*v \sim q = i_E(fw)i_E(fw)^* \sim i_E(fw)^*i_E(fw) = 1$
and $\Oo_E$ contains nonunitary isometries,
we deduce that $v^*v$ is an infinite projection in $A$.
\end{proof}

\section{Simplicity for lexicographic product systems}

For the lexicographic product systems of
Example~\ref{example:lexicographic},
we have the following partial converse to Theorem~\ref{theorem:main}.

\begin{theorem}\label{theorem:lexicographic}
Let $G$ be a countable abelian group,
let $P$ be a subsemigroup of $G$ which contains the identity,
let $d:P\to\NN^*$ be a semigroup homomorphism,
and let $E$ be the lexicographic product system determined by $d$.
Then $\Oo_E$ is simple if and only if $d$ is injective.
\end{theorem}

\begin{proof}
By Proposition~\ref{prop:nontrivial},
$E$ admits a nontrivial Cuntz representation,
so one direction follows from Theorem~\ref{theorem:main}.
For the converse, suppose $d$ is not injective.
Fix $s,t\in P$ such that $s\ne t$ and $d(s) = d(t)$.
We claim that $b := i_E(s,\delta_0) - i_E(t,\delta_0)$
generates a proper ideal $\Ii\triangleleft\Oo_E$,
so that $\Oo_E$ is not simple.
Let $S:E\to B(L^2(\TT))$ be the distinguished Cuntz representation
defined in Example~\ref{example:lexicographic}.
Since $S_*(b) = S(s,\delta_0) - S(t,\delta_0) = 0$
we have $b\in\ker S_*$,
and hence $\Ii$ is not all of $\Oo_E$.
To see that $\Ii$ is nonzero,
choose $\lambda\in\widehat G$ such that
$\lambda(s) \ne \lambda(t)$,
and define $T(r,x) := \lambda(r)S(r,x)$ for all $(r,x)\in E$.
Then $T$ is a Cuntz representation of $E$  such that
\[
T_*(b) = T(s,\delta_0) - T(t,\delta_0) = (\lambda(s) - \lambda(t))S(s,\delta_0) \ne 0,
\]
and we deduce that $b$, and hence $\Ii$, is nonzero.
\end{proof}

When $P = \NN^2$ we can say a bit more.
For every $m,n\in\NN^*$,
write $E(m,n)$ for lexicographic product system over $\NN^2$
determined by the homomorphism $(a,b)\in\NN^2 \mapsto m^an^b \in\NN^*$.
Note that $E(m,n) \cong E(n,m)$.

\begin{theorem}\label{theorem:N2}
If $m = 1$, then $\Oo_{E(m,n)} \cong  \Oo_n\otimes C(\TT)$.
If $m,n \ge 2$ and $\log_m n$ is irrational, then $\Oo_{E(m,n)}$
is simple.
If $m,n \ge 2$ and $\log_m n$ is rational, then
there exists a unique positive integer $l$
such that $m = l^a$ and $n = l^b$ 
for some relatively prime positive integers $a$ and $b$;
for this $l$ we have $\Oo_{E(m,n)} \cong \Oo_l\otimes C(\TT)$.
\end{theorem}

The proof of this Theorem is preceded by two Propositions.
The first deals with constructing representations of
lexicographic product systems over $\NN^k$,
allowing for $k = \infty$ by defining
$\NN^\infty := \bigoplus_{i=1}^\infty \NN$.
Let $\{e_a: 1 \le  a \le  k\}$ be the canonical basis for $\NN^k$.

\begin{prop}\label{prop:construct reps}
Let $d:\NN^k\to\NN^*$ be a semigroup homomorphism,
and let $B$ be a unital $C^*$\ndash algebra.
Suppose $\{U_{a,i}: 1 \le a \le k,\ 0 \le i \le d(e_a) - 1 \}$
is a set of isometries in $B$ which satisfies
\[
U_{a,i}^*U_{a,j} = 0
\qquad\text{whenever $i \ne j$, and}
\]
\begin{equation}\label{eq:commutation}
U_{a,i}U_{b,j} = U_{b,p}U_{a,q}
\qquad\text{whenever $id(e_b) + j = pd(e_a) + q$.}
\end{equation}
Then there is a unique representation $\phi:E(d) \to B$
such that $\phi(a,\delta_i) = U_{a,i}$,
and $\phi$ is a Cuntz representation if
\begin{equation}\label{eq:cuntz}
\sum_{i=0}^{d(e_a) - 1} U_{a,i}U_{a,i}^* = 1
\qquad\text{for every $a$.}
\end{equation}
\end{prop}

\begin{proof}
Fix $s\in\NN^k$ and $m\in\{0,\dots, d(s) - 1\}$.
Our first goal is to define $\phi(s,\delta_m)$.
In any expression of $s$ as an ordered sum
$e_{a_1} + \dotsb + e_{a_l}$ of basis elements
we have $l = \lVert s \rVert_1$.
For each of these finitely-many ordered $l$\ndash tuples
$(a_1, \dots, a_l)$,
it is easy to see that
there is a unique way of factoring $(s,\delta_m)$
as a product $(e_{a_1}, \delta_{i_1})\dotsb  (e_{a_l}, \delta_{i_l})$.
These factorizations can be obtained from one another
by using the commutation relations
\[
(e_a,\delta_i)(e_b, \delta_j) = (e_b, \delta_p)(e_a, \delta_q)
\qquad\text{whenever $id(e_b) + j = pd(e_a) + q$}
\]
on adjacent factors a finite number of times,
and hence  \eqref{eq:commutation} implies that the corresponding product
$U_{a_1,i_1}\dotsb U_{a_l,i_l}$ 
is independent of the factorization;
we define $\phi(s,\delta_m)$ to be this common element of $B$.
It is obvious that $\phi(s,\delta_m)\phi(t,\delta_n) = \phi((s,\delta_m)(t,\delta_n))$,
and defining
\[
\phi(s,x) := \sum_{i=0}^{d(s) - 1} (x \mid \delta_i) \phi(s, \delta_i)
\qquad\text{for $(s,x)\in E(d)$}
\]
gives the desired representation $\phi$.
If \eqref{eq:cuntz} holds, then
\[
\sum_{m=0}^{d(s) - 1} \phi(s,\delta_m)\phi(s,\delta_m)^*
= \sum_{i_1 = 0}^{d(e_{a_1}) - 1} \dotsb
  \sum_{i_l = 0}^{d(e_{a_l}) - 1}
  U_{a_1,i_1}\dotsb U_{a_l,i_l}U_{a_l,i_l}^*\dotsb U_{a_1,i_1}^*
= 1,
\]
so $\phi$ is a Cuntz representation.
\end{proof}

\begin{prop}\label{prop:factor}
$\Oo_{E(m,n)} \cong \Oo_{E(m, mn)}$.
\end{prop}

\begin{proof}
Let $E = E(m,n)$ and $F = E(m,mn)$.
Since $\dim F_{(a,b)} = m^a(mn)^b = m^{a+b}n^b = \dim E_{(a+b,b)}$,
we can define $\psi:F\to\Oo_E$ by
\[
\psi((a,b),x) = i_E((a+b,b),x)
\qquad\text{for all $((a,b),x)\in F$.}
\]
It is easy to see that $\psi$ is a Cuntz representation of $F$,
so there is a homomorphism $\psi_*:\Oo_F \to \Oo_E$
such that $\psi_*\circ i_F = \psi$.
Note that for any $(a,b)\in\NN^2$ and $j\in\{0,\dots, a^mb^n-1\}$
we have
\begin{align*}
i_E((a,b),\delta_j)
& = i_E((b,0),\delta_i)^*i_E((a+b,b),\delta_{im^an^b + j}) \\
& = \psi((b,0),\delta_i)^*\psi((a,b),\delta_{im^an^b + j});
\end{align*}
since elements of the form $i_E((a,b),\delta_j)$
generate $\Oo_E$, $\psi_*$ is surjective.

To see that $\psi_*$ is injective, we construct its inverse.
First define
\[
V_{1,i} := i_F((1,0),\delta_i) \in \Oo_F
\qquad\text{for $0 \le  i \le  m-1$, and}
\]
\[
V_{2,k} := i_F((0,1),\delta_k) \in \Oo_F
\qquad\text{for $0 \le  k \le  mn-1$.}
\]
Note that the $V_{1,i}$'s and $V_{2,k}$'s
are isometries which generate $\Oo_F$,
that $\sum V_{1,i}V_{1,i}^* = 1 = \sum V_{2,k}V_{2,k}^*$,
and that
\[
V_{1,i}V_{2,k} = V_{2,k'}V_{1,i'}
\qquad\text{whenever $imn + k =  k'm +  i'$.}
\]
Now define
$U_{1,i} := V_{1,i}$ for $0 \le i \le m-1$ and
\[
U_{2,j} := \sum_{l=0}^{m-1} V_{2,jm+l}V_{1,l}^*
\qquad\text{for $0 \le j \le n-1$.}
\]
It routine to check that the $U_{2,j}$'s are isometries whose  range projections
sum to the identity.
Suppose $in + j = pm + q$, where
$0 \le i,q \le  m-1$ and $0 \le j,p \le n-1$.  Then
\begin{align*}
U_{1,i}U_{2,j}
& = \sum_{l=0}^{m-1} V_{1,i}V_{2,jm + l}V_{1,l}^*
  = \sum_{l=0}^{m-1} V_{2,in + j}V_{1,l}V_{1,l}^* \\
& = V_{2,in + m} = V_{2,pm + q}
  = \sum_{l=0}^{m-1} V_{2,pm + l}V_{1,l}^*V_{1,q}
  = U_{2,p}U_{1,q},
\end{align*}
so by Proposition~\ref{prop:construct reps} there is a
Cuntz representation $\phi:E\to\Oo_F$ such that
$\phi((1,0),\delta_i) = U_{1,i}$ and $\phi((0,1),\delta_j) = U_{2,j}$.
We check that $\phi_*\circ \psi_*$ is the identity on $\Oo_F$,
from which it follows that $\psi_*$ is injective, and hence an isomorphism:
if $0 \le i \le m-1$, then
\[
\phi_*\circ\psi_*(V_{1,i})
= \phi_*(\psi((1,0),\delta_i))
= \phi_*(i_E((1,0),\delta_i))
= U_{1,i} = V_{1,i};
\]
if also $0 \le j \le n-1$, then
\begin{align*}
\phi_*\circ\psi_*(V_{2,in + j})
& = \phi_*(\psi((0,1),\delta_{in + j})) \\
& = \phi_*(i_E((1,1),\delta_{in + j}))
= U_{1,i}U_{2,j} = V_{2,in + j}.
\end{align*}
\end{proof}

\begin{proof}[Proof of Theorem~\ref{theorem:N2}]
By Proposition~\ref{prop:construct reps},
$\Oo_{E(1,n)}$ is the universal $C^*$\ndash al\-ge\-bra for
collections $\{U_{1,0}\} \cup \{U_{2,0}, \dots, U_{2,n-1}\}$
of isometries satisfying
\[
U_{1,0}U_{1,0}^* = 1,
\quad
\sum_{i=0}^{n-1} U_{2,i}U_{2,i}^* = 1,
\quad\text{and}\quad
U_{1,0}U_{2,i} = U_{2,i}U_{1,0};
\]
i.e., $\Oo_{E(1,n)} \cong \Oo_n\otimes C(\TT)$.
Suppose $m,n\ge 2$.  If $\log_m n$ is irrational,
then $d: (a,b)\in\NN^2 \mapsto m^an^b$ is injective,
and $\Oo_{E(m,n)}$ is simple by Theorem~\ref{theorem:lexicographic}.
The existence and uniqueness of $l$ is elementary when $\log_m n$ is rational,
and repeated applications of Proposition~\ref{prop:factor}
give $\Oo_{E(m,n)} \cong \Oo_{E(l,1)}$.
\end{proof}

\section{Nuclearity}
Following \cite{dinhjfa,laca}, we call a cross section
$u:s\in P\mapsto u_s\in E_s\setminus\{0\}$
a {\em twisted unit\/} of $E$ if $u_su_t \in \CC u_{st}$ for every $s,t\in P$.

\begin{theorem}\label{theorem:nuclearity}
Let $G$ be a countable abelian group,
let $P$ be a subsemigroup of $G$ which contains the identity,
and let $E$ be a product system over $P$ of finite-dimensional Hilbert spaces.
If $E$ admits a twisted unit, then $\Oo_E$ is nuclear.
\end{theorem}

Our proof is modelled on the one given by Laca in \cite[Section~3]{laca}:
we realize $\Oo_E$ as a twisted semigroup crossed product of an AF algebra
by $P$, and deduce nuclearity from a theorem of Murphy
\cite[Theorem~3.1]{murphy}.
We begin by recalling the definition of a twisted semigroup crossed product.
Let $P$ be as in the theorem,
let $\beta$ be an action of $P$ as endomorphisms
of a unital $C^*$\ndash algebra $A$,
and let $\omega$ be a multiplier of $P$.
We call $(A,P,\beta,\omega)$ a {\em twisted semigroup dynamical system\/}.

Suppose $B$ is a unital $C^*$\ndash algebra.
An {\em isometric $\omega$\ndash representation\/}
of $P$ in $B$ is a map $V:P\to B$
such that each $V_s$ is an isometry and $V_sV_t = \omega(s,t)V_{st}$
for all $s,t\in P$.
A {\em covariant representation\/} of $(A,P,\beta,\omega)$ in $B$
is a pair $(\pi,V)$ consisting of a unital homomorphism $\pi:A\to B$
and an isometric $\omega$\ndash representation $V:P\to B$ such that
\[
\pi(\beta_s(a)) = V_s\pi(a)V_s^*
\qquad\text{for all $s\in P$ and $a\in A$.}
\]
A {\em crossed product\/} for $(A,P,\beta,\omega)$ is a triple
$(C,i_A,i_P)$ consisting of a unital $C^*$\ndash algebra $C$
and a covariant representation $(i_A,i_P)$ of $(A,P,\beta,\omega)$ in $C$
such that

(a) for every covariant representation $(\pi,V)$ of $(A,P,\beta,\omega)$
in a unital $C^*$\ndash algebra $B$,
there is a homomorphism $\pi\times V:C\to B$ such that
$(\pi\times V)\circ i_A = \pi$ and $(\pi\times V)\circ i_P = V$; and

(b) $C$ is generated as a $C^*$\ndash algebra by $i_A(A)\cup i_P(P)$.

\noindent The triple $(C,i_A,i_P)$ is unique up to canonical isomorphism.

\begin{proof}[Proof of Theorem~\ref{theorem:nuclearity}]
Let $u$ be a twisted unit for $E$.
By replacing $u_s$ with $\lVert u_s \rVert^{-1}u_s$,
we can assume that each $u_s$ is a unit vector.
Then $u_su_t = \omega(s,t)u_{st}$ determines a multiplier $\omega$ of $P$.

As in Section~\ref{section:abelian},
let $\Ff_E$ be the inductive limit $\varinjlim\Kk(E_s)$
under the embeddings $S\in\Kk(E_s) \mapsto S\otimes 1^t\in\Kk(E_{st})$,
and let $\iota_s$ be the canonical embedding of $\Kk(E_s)$ in $\Ff_E$.
Since tensoring on the left
by the rank-one projection $\rankone{u_r}{u_r}$
commutes with tensoring on the right by the identity, 
for each $r\in P$
there is an endomorphism $\beta_r$ of $\Ff_E$ which satisfies
\[
\beta_r(\iota_s(S)) = \iota_{rs}((\rankone{u_r}{u_r})\otimes S)
\qquad\text{for all $s\in P$ and $S\in\Kk(E_s)$.}
\]
Note that $\beta_r$ is injective.
Moreover, for  any $a\in \Ff_E$ we have
\begin{align*}
\beta_r\circ\beta_s(a)
& = (\rankone{u_r}{u_r}) \otimes (\rankone{u_s}{u_s}) \otimes a
  = (\rankone{u_ru_s}{u_ru_s}) \otimes a \\
& = (\rankone{\omega(r,s)u_{rs}}{\omega(r,s)u_{rs}}) \otimes a
  = (\rankone{u_{rs}}{u_{rs}}) \otimes a
  = \beta_{rs}(a),
\end{align*}
and hence $\beta:P\to\End \Ff_E$ is a semigroup homomorphism.

Let $i_\Ff:\Ff_E\to\Oo_E$ be the embedding of
Proposition~\ref{prop:FE},
and define $i_P:P\to\Oo_E$ by $i_P(s) := i_E(u_s)$.
We claim that $(\Oo_E,i_\Ff,i_P)$ is a crossed product for the twisted
semigroup dynamical system $(\Ff_E,P,\beta,\omega)$.
Each $i_P(s)$ is an isometry,
and
$i_P(s)i_P(t) = i_E(u_su_t) = \omega(s,t)i_E(u_{st}) = \omega(s,t)i_P(st)$,
so $i_P$ is an isometric $\omega$-representation of $P$ in $\Oo_E$.
If $x,y\in E_s$, then
\begin{align*}
i_\Ff(\beta_r(\iota_s(\rankone xy)))
& = i_\Ff(\iota_{rs}((\rankone{u_r}{u_r})\otimes(\rankone xy)))
  = i_\Ff(\iota_{rs}(\rankone{u_rx}{u_ry})) \\
& = i_E(u_rx)i_E(u_ry)^*
  = i_P(r)i_\Ff(\iota_s(\rankone xy))i_P(r)^*,
\end{align*}
and since elements of the form $\iota_s(\rankone xy)$ have dense
linear span in $\Ff_E$, we deduce that
$i_\Ff(\beta_r(a)) = i_P(r)i_\Ff(a)i_P(r)^*$
for all $r\in P$ and $a\in \Ff_E$.
Thus $(i_\Ff,i_P)$ is a covariant representation of $(\Ff_E,P,\beta,\omega)$
in $\Oo_E$.

We now verify condition (a) of a crossed product.
Suppose $(\pi,V)$ is a covariant representation of
$(\Ff_E,P,\beta,\omega)$ in a unital $C^*$\ndash algebra $B$.
Define $\phi:E\to B$ by
\[
\phi(x) := \pi(\iota_{p(x)}(\rankone x{u_{p(x)}}))V_{p(x)}
\qquad\text{for $x\in E$.}
\]
We claim that $\phi$ is a Cuntz representation of $E$ in $B$.
If $x,y\in E_s$, then
\begin{align*}
\phi(y)^*\phi(x)
& = V_s^*\pi(\iota_s(\rankone y{u_s}))^*\pi(\iota_s(\rankone x{u_s}))V_s \\
& = V_s^*\pi(\iota_s((\rankone {u_s}y)(\rankone x{u_s})))V_s \\
& = (x\mid y)V_s^*\pi(\iota_s(\rankone{u_s}{u_s}))V_s \\
& = (x\mid y)V_s^*\pi(\beta_s(1))V_s = (x\mid y)1,
\end{align*}
and
\begin{equation}\label{eq:xystar}
\begin{split}
\phi(x)\phi(y)^*
& = \pi(\iota_s(\rankone x{u_s}))V_sV_s^*\pi(\iota_s(\rankone y{u_s}))^* \\
& = \pi(\iota_s((\rankone x{u_s})(\rankone{u_s}{u_s})(\rankone {u_s}y)))
  = \pi(\iota_s(\rankone xy)).
\end{split}
\end{equation}
With $\Bb_s$ an orthonormal basis for $E_s$, \eqref{eq:xystar} gives
\[
\sum_{f\in \Bb_s} \phi(f)\phi(f)^*
= \sum_{f\in \Bb_s} \pi(\iota_s(\rankone ff))
= \pi(\iota_s(1)) = 1 = \phi(1).
\]
To see that $\phi$ is multiplicative,
suppose $x\in E_s$ and $y\in E_t$.
Then
\begin{align*}
\phi(x)\phi(y)
& = \pi(\iota_s(\rankone x{u_s}))V_s\pi(\iota_t(\rankone y{u_t}))V_t \\
& = \pi(\iota_s(\rankone x{u_s}))\pi(\beta_s(\iota_t(\rankone y{u_t})))V_sV_t \\
& = \pi(\iota_{st}((\rankone x{u_s})\otimes 1^t)
                   ((\rankone{u_s}{u_s})\otimes(\rankone y{u_t})))V_sV_t \\
& = \pi(\iota_{st}((\rankone  x{u_s})\otimes(\rankone y{u_t})))V_sV_t \\
& = \pi(\iota_{st}(\rankone{xy}{u_su_t}))V_sV_t \\
& = \pi(\iota_{st}(\rankone{xy}{u_{st}}))\overline{\omega(s,t)}V_sV_t \\
& = \pi(\iota_{st}(\rankone{xy}{u_{st}}))V_{st}
  = \phi(xy).
\end{align*}
Thus $\phi$ is a Cuntz  representation of $E$.

Let $\pi\times V:\Oo_E\to B$ be the integrated form of $\phi$;
that is, $\pi\times V$ satisfies
$(\pi\times V)\circ i_E = \phi$.
We claim that $(\pi\times V)\circ i_\Ff = \pi$
and $(\pi\times V)\circ i_P = V$, giving condition (a) of a crossed product.
If $x,y\in E_s$, then by \eqref{eq:xystar} we have
\[
(\pi\times V)\circ i_\Ff(\iota_s(\rankone xy))
= \pi\times V(i_E(x)i_E(y)^*)
= \phi(x)\phi(y)^*
= \pi(\iota_s(\rankone xy)),
\]
and since  elements of the form $\iota_s(\rankone xy)$ have dense linear span
in $\Ff_E$, we deduce that $(\pi\times V)\circ i_\Ff = \pi$.
For the second part of the claim, we calculate
\begin{align*}
(\pi\times V)\circ i_P(s)
& = \pi\times V(i_E(u_s))
  = \phi(u_s) \\
& = \pi(\iota_s(\rankone{u_s}{u_s}))V_s
  = \pi(\beta_s(1))V_s = V_sV_s^*V_s = V_s.
\end{align*}

To verify condition (b) of a crossed product, suppose $x\in E_s$.
Then
\[
i_E(x) = i_E(x)i_E(u_s)^*i_E(u_s)
= i_\Ff(\iota_s(\rankone x{u_s}))i_P(s) \in C^*(i_\Ff(\Ff_E)\cup i_P(P)),
\]
and since $\Oo_E$ is generated as a $C^*$\ndash algebra by $i_E(E)$,
this shows that it is also generated by $i_\Ff(\Ff_E)\cup i_P(P)$, as required.

We have shown that $(\Oo_E,i_\Ff,i_P)$ is a crossed product for $(\Ff_E,P,\beta,\omega)$.
Since $\Oo_E$ is the twisted semigroup crossed  product
of a nuclear $C^*$\ndash algebra by an abelian semigroup of injective endomorphisms,
it is nuclear by \cite[Theorem~3.1]{murphy}.
\end{proof}

\begin{cor}\label{cor:nuclearity}
If $E$ is a twisted lexicographic product system over a subsemigroup
of a countable abelian group, then $\Oo_E$ is nuclear.
\end{cor}

\begin{proof} Since $u_s := (s,\delta_0)$ is a twisted unit for $E$,
Theorem~\ref{theorem:nuclearity} applies.
\end{proof}

\begin{remark} Suppose $P$ is a subsemigroup of a countable abelian group,
$d:P\to\NN^*$ is an injective homomorphism, and $\omega$ is a multiplier
of $P$.  Let $E = E(d)^\omega$ be the lexicographic product system $E(d)$
twisted by $\omega$.  Then $\Oo_E$ is a separable, unital
$C^*$\ndash algebra which is simple  and purely infinite
(Theorem~\ref{theorem:main}) and nuclear (Corollary~\ref{cor:nuclearity}).
If $P = \NN^k$ and $\omega$ is trivial (i.e., $E = E(d)$),
then by the proof of Theorem~\ref{theorem:nuclearity}
$\Oo_E$ is the crossed product of an AF algebra by $\NN^k$,
hence a full corner of the crossed product of an AF algebra by $\ZZ^k$,
and hence belongs to the bootstrap class of algebras which satisfy
the Universal Coefficient Theorem \cite[Theorem~23.1.1]{blackadar}.
\end{remark}


\begin{thebibliography}{20}

\bibitem{arv} W.\ Arveson, {\em Continuous analogues of Fock space\/},
Memoirs Amer.\ Math.\ Soc.\ {\bf 80} (1989), No.~409.

\bibitem{bprs} T.\ Bates, D.\ Pask, I.\ Raeburn and W.\ Szyma\'nski,
{\em The $C^*$-algebras of row-finite graphs}, preprint.

\bibitem{blackadar} B.\ Blackadar, {\em K-theory for Operator Algebras,
Second Edition\/}, MSRI Publications~5, Cambridge University Press,
Cambridge, 1998.

\bibitem{cun} J.\ Cuntz, {\em  Simple $C^*$-algebras generated by isometries\/},
  Comm.\ Math.\ Phys.\ {\bf 57} (1977), 173--185.

\bibitem{dinhjfa} H.\ T.\ Dinh, {\em Discrete product systems and
their $C^*$-algebras\/}, J.\ Funct.\ Anal.\ {\bf 102} 
(1991), 1--34.

\bibitem{dinhjot} H.\ T.\ Dinh, {\em On generalized Cuntz $C^*$-algebras\/},
J.\ Operator Theory {\bf 30} (1993), 123--135.

\bibitem{fowler} N.\ J.\ Fowler, {\em Compactly-aligned discrete product systems,
and generalizations of $\mathcal O_\infty$\/}, International J. Math., to appear.

\bibitem{fowrae} N.\ J.\ Fowler and I.\ Raeburn, {\em Discrete product systems
and twisted crossed products by semigroups\/}, J.\ Funct.\ Anal.\ {\bf 155}
(1998), 171--204.

\bibitem{fowrae2} N.\ J.\ Fowler and I.\ Raeburn,
{\em The Toeplitz algebra of a Hilbert bimodule\/},
Indiana Univ.\ Math.\ J., to appear.

\bibitem{laca} M.\ Laca, {\em Discrete product systems with twisted units\/},
Bull.\ Austral.\ Math.\ Soc.\ {\bf 52} (1995), 317--326.

\bibitem{laca2} M.\ Laca, {\em From endomorphisms to automorphisms and back:
dilations and full corners\/}, preprint.

\bibitem{murphy} G.\ J.\ Murphy, {\em Crossed products of $C^*$-algebras
by endomorphisms\/}, Int.\ Eq.\ and Operator Theory {\bf 24} (1996), 298--319.

\end{thebibliography}
\end{document}